\documentclass{article}

\usepackage{amsmath,amssymb,amsfonts}

\newtheorem{theorem}{Theorem}

\newcommand\D{{\cal D}}
\newcommand\R{\mathbb R}
\newcommand\C{\mathbb C}
\newcommand\nil{{\mathrm{Nil}}}
\newcommand\SU{{\mathrm{SU}}}
\newcommand\const{{\mathrm{const}}}
\newcommand\tr{{\mathrm{tr}}}
\newcommand\ovl[1]{{}\,\overline{\!#1}\,}
\newcommand\ovp[1]{{}\,\overline{{\!#1}}\,}
\newcommand\ovw{{}\,\overline{\hskip-1mm\widetilde{A}}}
\newcommand\wta{{}\,\widetilde{\!A}}
\newcommand\wh{{}\,\widehat{h}}
\newcommand\dv{\kern.3mm{\rm{:}}\ }

\begin{document}

\title{On Constant Mean Curvature Surfaces in the Heisenberg Group
%\thanks{This work is supported by Russian %Foundation for Basic Research (project %number~09-01-00598),
%State Maintenance Program for the Leading %Scientific Schools of the
%Russian Federation (grant ~Nsh-7256.2010.1).}
}
\author{Berdinsky~D.
\thanks{
%Sobolev Institute of Mathematics, 630090 %Novosibirsk, Russia; 
email: berdinsky@gmail.com.}  }

\date{}
\maketitle

\begin{quote}

\noindent{\bf Abstract}. 
%%%%%%%%%%%%%%%%%%%%%%%%%%%
% old version 
%%%%%%%%%%%%%%%%%%%%%%%%%%%%
%This paper is devoted to the %theory of surfaces of %constant mean curvature in %the three--dimensional %Heisenberg group. We prove %that each surface of such a %kind locally corresponds
%to some solution of the %system of a sine--Gordon %type equation and a first %order partial differential %equation. 
%In this work 
We study constant mean curvature surfaces in the three--dimensional Heisenberg group. We prove that a % each CMC 
constant mean curvature surface in a neighborhood of non--umbilic point is described by some solution of 
a sinh--Gordon equation subject to a first order differential constraint.

\noindent{\bf Keywords}: Heisenberg group, surfaces of constant mean curvature 

\end{quote}

\section{Introduction}

In this paper we study constant mean curvature
surfaces in the three--dimensional Heisenberg group~$\nil$ with the Thurston metric using the Weierstrass representation %technique for 
of surfaces in three--dimensional Lie groups. 
The Weierstrass representation of surfaces 
in $\SU(2)$  is 
described in \cite{T}.
For surfaces in the three--dimensional noncompact Lie groups with the Thurston metrics, including the Heisenberg group $\nil$, 
the Weierstrass representation  is 
described in~\cite{BT05}.
%This representation was proposed
%in~\cite{T}, where, as an example, it is %implemented for surfaces
%in~$\SU(2)$. For surfaces in the %three--dimensional noncompact Lie
%groups with the Thurston metrics 
%the Weierstrass representation technique is %described in~\cite{BT05}.
%The Abresch--Rosenberg quadratic differential %is holomorphic on constant mean curvature %surfaces in $\nil$~\cite{AR05}. 

The {\it Abresch--Rosenberg quadratic differential} 
for surface in $\nil$ is a quadratic 
differential which is holomorphic %exactly 
on constant mean curvature surfaces~\cite{AR05}. 
We say that a point on a surface 
in $\nil$ is {\it umbilic} %umbilical
if the Abresch--Rosenberg quadratic vanishes 
at this point.   
%%%%%%%%%%%%%%%%%%%%%%%%
%We define an umbilic
%point of a surface as a point where the %Abresch--Rosenberg quadratic differential %vanishes.
%%%%%%%%%%%%%%%%%%%%%%%%
%In the sequel we will consider constant mean %curvature surfaces only in
%neighborhoods of the non--umbilical points. 
%%%%%%%%%%%%%%%%%%%%%%%%
In this paper 
we consider constant mean  
curvature surfaces in neighborhoods of 
non--umbilic points. 
%%%%%%%%%%%%%%%%%%%%%%%%
%If a surface is compact and
%all its points are umbilical then it is a %sphere of constant mean curvature. 
%%%%%%%%%%%%%%%%%%%%%%%%
A closed surface in $\nil$ for which all 
its points are umbilic is a constant 
mean curvature sphere 
(for given mean curvature $H \neq 0$ it is unique up to isometry).   
%%%%%%%%%%%%%%%%%%%%%%%%
%Herewith there are no umbililical points on a %torus of constant mean curvature. 
%%%%%%%%%%%%%%%%%%%%%%%%
A torus of a constant mean curvature in
$\nil$ (if it exists) does not have umbilic points. The original motivation of this paper
stemmed from the problem of finding immersions 
of constant mean curvature tori in $\nil$.   

%%%%%%%%%%%%%%%%%%%%%%%%%%%%%%%%%%%%%%%%%%%
% The important distinction (already noted %in~\cite{T}) of the Weierstrass
% representation for surfaces in the %non--commutative Lie groups from that 
% for surfaces in~$\R^3$ is the fact that if a
% spinor~$\psi$ satisfies ${\cal D}\psi=0$ then %the spinor~$r\psi$, $r\in\R, r\ne\pm 1$
% does not necessary correspond to some surface.
% In other words, in general case, the solution %of some Dirac equation does not correspond to %some immersion of a surface (in
% contrast to the Weierstrass representation %for surfaces in~$\R^3$).
%%%%%%%%%%%%%%%%%%%%%%%%%%%%%%%%%%%%%%%%%%%
%In~\S\,2 we show that for a neighborhood of a %non--umbilical point
%of a constant mean curvature surface a %sine--Gordon type
%equation is satisfied for $v=\log U$ where %$U$~ is the potential of
%the Dirac operator corresponding to the %Weirstrass representation \big(it was
%proved in~2007 and briefly mentioned %in~\cite{T08}\big).
%%%%%%%%%%%%%%%%%%%%%%%%%%%%%%%%%%%%%%%%%%%

%We show that for a neighborhood of a %non--umbilical point
%of a constant mean curvature surface a %sine--Gordon type
%equation is satisfied for $v=\log U$ where %$U$~ is the potential of
%the Dirac operator corresponding to the %Weirstrass representation \big(it was
%proved in~2007 and briefly mentioned %in~\cite{T08}\big).

\emph{Contributions of the paper.}
We will show that a surface of 
constant mean curvature $H \neq 0$ 
in $\nil$
in a neighborhood of a non--umbilic 
point corresponds to a solution  of 
a sinh--Gordon equation of the form:
\begin{equation}
\label{sinh-Gordon}
\Delta v+ 8 \sinh 2v=0,
\end{equation}
subject to a first 
order differential constraint\footnote{To derive this first order differential constraint we make a couple of assumptions 
about the geometry of the surface, see Theorem \ref{theocmcsyst}.}: 
%(we make few additional assumptions for $v$,
%see Theorem \ref{theocmcsyst}): 
%(referred to as the {\it reality condition})
\begin{equation} 
\label{first_order_constr}   
   \frac{(\operatorname{Im} (v))_x^2}
   {\cosh^2 \operatorname{Re}(v)} + 
   \frac{(\operatorname{Im} (v))_y^2}
   {\sinh^2 \operatorname{Re}(v)} =
   8
   \left(\cos (2\operatorname{Im}(v))-\operatorname{Re}
   \left(\dfrac{2H+i}{2H-i}\right) \right),
\end{equation}	
where $v = \log U$ and $U=U(z)$ is the complex--valued potential 
of the Dirac operator emerging from 
the Weierstrass representation of 
the surface with respect a conformal parameter $z= x+ iy$\,\footnote{This result was proved 
in 2007 and it was mentioned in \cite{B09,T08}.}.
Moreover, a minimal surface $H=0$ in 
$\nil$ in a neighborhood of a 
non--umbilic point corresponds 
to a solution of the equation
\eqref{sinh-Gordon} satisfying 
the constraint $\operatorname{Re}(e^v)=0$:
\begin{equation} 
\label{system_H_equals_zero}  
  \Delta v+ 8 \sinh 2v=0,\,\,\,\, 
  \operatorname{Re}(e^v)=0.	
\end{equation}	 
Similarly to the three--dimensional Euclidean
space, the sinh--Gordon equation in the 
system \eqref{sinh-Gordon}--\eqref{first_order_constr}
and the system \eqref{system_H_equals_zero}
enables us to use methods of soliton theory 
for constructing immersions of constant 
mean curvature surfaces in $\nil$. See, e.g., 
the use of Lamb ansatz for constructing immersions of constant mean curvature tori in the three--dimensional Euclidean space \cite{Abresch87}.  

{\it The structure of the paper.}
The rest of the paper is organized as follows. 
In Section \ref{Gauss-Weingarten_section} 
for a constant mean curvature surface  
in $\nil$ we derive the Gauss--Weingarten equations and its compatibility condition. 
In Section \ref{section_reality_condition}
we derive the additional constraint for the potential of a Dirac operator 
that is necessary for the existence of a constant mean curvature surface in $\nil$ corresponding to this potential. Theorem 
\ref{theocmcsyst} and Theorem \ref{minimal_surf_thm} summarize 
the obtained results for nonzero 
constant mean curvature surfaces and
minimal surfaces, respectively.

%in a neighborhood of a non--umbilic point
%on  a sine--Gordon type equation is satisfied %for $v=\log U$ where $U$~ is the potential of
%the Dirac operator corresponding to the %Weirstrass representation 
%\big(it was proved in~2007 and briefly %mentioned in~\cite{T08}\big).

%In~\S\,3 we study the question of searching a %solution to a sine--Gordon type equation such %that, for the corresponding Dirac
%operator~${\cal D}$, there exists a solution %to the equation  ${\cal D}\psi=0$
%corresponding to some constant mean curvature %surface. We deduce the equation (called the %reality condition) with the desired
%property of correspondence to some constant %mean curvature surface.

%The equations obtained  enable us to use the 
%Lamb anzats technique
%\big(see~\cite{Abresch87}\big) or other %methods of the soliton theory
%to construct an example of immersing of a %constant mean curvature
%torus in~$\nil$.

%%%%%%%%%%%%%%%%%%%%%%%%%%%%%%%%%%%%%%%%%%%%%% 

%In~\S\,3 we study the question of searching a %solution to a sine--Gordon type equation such %that, for the corresponding Dirac
%operator~${\cal D}$, there exists a solution %to the equation  ${\cal D}\psi=0$
%corresponding to some constant mean curvature %surface. We deduce the equation (called the %reality condition) with the desired
%property of correspondence to some constant %mean curvature surface.

\section{The Gauss--Weingarten Equations and its Compatibility Condition}
\label{Gauss-Weingarten_section}

In this section we derive 
the Gauss--Weingarten equations and its 
compatibility condition 
for a constant mean curvature surface 
in $\nil$. We express 
these equations 
in terms of  
{\it generating spinors} of the surface, 
the mean curvature, the 
Abresch--Rosenberg quadratic differential and a logarithm of the potential of the Dirac operator emerging from the Weierstrass 
representation  of the surface; see the 
equations \eqref{gauss_weingraten_eq} 
and \eqref{compcond}.

We represent the three--dimensional Heisenberg group~$\nil$ with the Thurston metric~\cite{Scott} as the group of  upper
triangular matrices of the form:
$$
\begin{pmatrix}
1 & x_1 & x_3\\
0 & 1 & x_2\\
0 & 0 & 1
\end{pmatrix},
\quad x_1,x_2,x_3\in\R,
$$
with the standard multiplication rule  and the left--invariant metric:
$$
ds^2=dx_1^2+dx_2^2+(dx_3-x_1dx_2)^2.
$$
Every surface $\mathcal{S}$ in $\nil$ can be locally  represented as a conformal mapping 
$$f:\Sigma\rightarrow\nil$$ 
from an open domain $\Sigma\subset\C$ 
into $\nil$. The coordinate $z \in \Sigma$
is a conformal parameter on $\mathcal{S}$.   
%defining a conformal immersion of a
%domain $\Sigma\subset\C$ into~$\nil$. The %conformal parameter on a
%surface is a coordinate $z\in\Sigma$.
Now we briefly recall the Weierstrass representation of surfaces in $\nil$ 
described in \cite{BT05}. 

Let $e_1,e_2,e_3$ be the  orthonormal basis of the Lie algebra for $\nil$ defined at
the identity of $\nil$ as:
$$e_1=\dfrac{\partial}{\partial x_1},\,\,\,\,
 e_2=\dfrac{\partial}{\partial x_2},\,\,\,\,
 e_3=\dfrac{\partial}{\partial x_3}.$$
The surface $S$ in $\nil$ can be
obtained 
by solving the equation:
\begin{equation}
\label{decomposition_f_z}	
f^{-1} \frac{\partial f}
{\partial z}=
\frac{i}{2}\big(\ovl{\psi}_2^2+\psi_1^2
\big)e_1+
\frac{1}{2}\big(\ovl{\psi}_2^2-\psi_1^2
\big)e_2+\psi_1\ovl{\psi_2} e_3,
\end{equation}
where the complex--valued functions $\psi_1(z)$, $\psi_2(z)$ satisfy the Dirac equation:
\begin{equation}
\gathered\label{nil-dirac}                        %(1)
\D\psi= \left[\begin{pmatrix}
         0 & \partial\\
        -\ovp{\partial}& 0
      \end{pmatrix}+
      \begin{pmatrix}
         U & 0\\
         0 & U
\end{pmatrix}
\right]
\psi=0,\\
\psi= \left(\begin{smallmatrix}
	\psi_1\\
	\psi_2
\end{smallmatrix}
\right),\,\,\,\,
U=\frac{H}{2} \big(|\psi_1|^2+|\psi_2|^2 \big)+\frac{i}{4}
\big(|\psi_2|^2-|\psi_1|^2 \big)
\endgathered
\end{equation}
and $H$ is the mean curvature of $\mathcal{S}$.
We will refer to 
$\psi= \left(\begin{smallmatrix}
	\psi_1\\
	\psi_2
\end{smallmatrix}
\right)$ as a  generating spinor
of the surface $\mathcal{S}$.

Let $\wta$ be the following 
complex--valued function of $z$ defined 
in terms of $\psi_1, \psi_2$ and $H$: 
\begin{equation}\label{HolDif}                    %(2)
	\wta= \ovl{\psi_2}\partial\psi_1-\psi_1
	\partial
	\ovl{\psi_2}+\frac{2Hi}{2H+i}\psi_1^2 \ovl{\psi}_2^2,
\end{equation}
The quadratic differential $\wta dz^2$ on a surface $\mathcal{S}$ is holomorphic if and only if the  surface $\mathcal{S}$ is of constant mean curvature \cite{BT05}. 
%This differential  $\wta dz^2$ is 
%the Abresch--Rosenberg quadratic %differential~\cite{AR05}. 
%The equation \eqref{HolDif} gives the %representation of the Abresch--Rosenberg 
%quadratic differential in terms of~$\psi_1$ %and~$\psi_2$.

From now on we assume that the mean curvature
$H$ of $\mathcal{S}$ is constant: $H=\const$. 
%It implies that
Therefore, the quadratic differential  $\wta dz^2$ is holomorphic. 
%Put $e^v:=U$. 
We denote by $v$ a logarithm of $U$: 
$$
  e^v = U. 
$$
%Using~\eqref{nil-dirac}, we represent the %derivative
It follows from \eqref{nil-dirac} 
that the derivative 
$\dfrac{\partial}{\partial z}U = v_z e^v$ has the form:
\begin{equation}
\label{U over z}                  %(3)
\frac{\partial}
     {\partial z}U=v_z e^v=
\frac14(2H+i)\psi_2\partial\ovl{\psi_2}+
\frac14(2H-i)\ovl{\psi_1}\partial\psi_1-
\frac{iH}2\psi_1\ovl{\psi_2}|\psi_2|^2.
\end{equation}
%Relations~\eqref{HolDif} and~\eqref{U over z} %give us the following form of the derivative
From the identities \eqref{HolDif} and 
\eqref{U over z} we obtain that: 
\begin{equation}\label{dpsi1}                     %(4)
\partial\psi_1=v_z\psi_1+
\frac14(2H+i)\wta e^{-v}\psi_2,
\end{equation}
%It is easy to verify that
%$$
and 
\begin{equation}\label{dpsi2}                     %(5)
\ovp{\partial}\psi_2=-
\frac14(2H-i)e^{-v}\ovw\psi_1+v_{\ovl{z}}\psi_2.
\end{equation}
Let $$B=\frac14(2H+i)\wta.$$ 
We refer to $Bdz^2$ as the 
Abresch--Rosenberg quadratic differential.
Similarly to $\wta dz^2$, $Bdz^2$ is 
holomorphic %exactly 
on constant mean curvature surfaces in $\nil$.

The equations~\eqref{nil-dirac},
\eqref{dpsi1} and~\eqref{dpsi2} enable us to represent the Gauss--Weingarten equations as follows:
\begin{equation}
\label{gauss_weingraten_eq}
\partial
\begin{pmatrix}
\psi_1\\
\psi_2
\end{pmatrix}
=\begin{pmatrix}
        v_z & Be^{-v}\\
       -e^v & 0
\end{pmatrix}
\begin{pmatrix}
       \psi_1\\
       \psi_2
\end{pmatrix},
\quad \ovp{\partial}
\begin{pmatrix}
       \psi_1\\
       \psi_2
\end{pmatrix}
=\begin{pmatrix}
       0 & e^v\\
       -\ovl{B}e^{-v} & v_{\ovl{z}}
\end{pmatrix}
\begin{pmatrix}
       \psi_1\\
       \psi_2
\end{pmatrix}.
\end{equation}
It %is easy to verify 
can be verified that 
that holomorphicity of the 
quadratic differential 
$Bdz^2$ and the compatibility condition 
of the system \eqref{gauss_weingraten_eq} 
yield the identity:
\begin{equation}\label{compcond}                  %(6)
v_{z\ovl{z}}+e^{2v}-|B|^2e^{-2v}=0.
\end{equation}

\section{The Reality Condition}               
\label{section_reality_condition}

 We are interested in the solutions $v$ of the equation \eqref{compcond} 
 for which there
 exists a generating spinor 
$\psi= \left(\begin{smallmatrix}
\psi_1\\
\psi_2
\end{smallmatrix}
\right)$ 
that satisfies the Gauss--Weingarten equations \eqref{gauss_weingraten_eq} together %and 
with the following identity:  
\begin{equation}\label{cond1}                     %(7)
e^{v}= \frac{H}2\left(|\psi_1|^2+|\psi_2|^2
         \right)+
\frac{i}4\left(|\psi_2|^2-|\psi_1|^2
         \right).
\end{equation}
Note that if the 
Gauss--Weingarten equations 
\eqref{gauss_weingraten_eq} and 
the identity \eqref{cond1} hold, then the identity \eqref{HolDif} holds as well.  
In this section we will derive a 
necessary condition for the 
existence of the generating 
spinor $\psi$ that satisfies
the Gauss--Weingarten equations 
\eqref{gauss_weingraten_eq} and the identity 
\eqref{cond1}.  

%Now we represent 
Let us represent 
the identity \eqref{cond1} as follows:
\begin{equation}\label{v_cond}                    %(8)
1=\left(\begin{array}{cc}
        \overline{\psi}_1&\overline{\psi}_2
        \end{array}
  \right)
  \left(\begin{array}{cc}
        \frac14(2H-i)e^{-v}&0\\
         0&\frac{1}{4}(2H+i) e^{-v}
        \end{array}
  \right)
  \left(\begin{array}{c}
        \psi_1\\
        \psi_2
        \end{array}
  \right).
\end{equation}
%Denote by~
Let $\wh$ be the following matrix:
$$
\wh= \left(\begin{array}{cc}
      \frac{1}{4}(2H-i) e^{-v} & 0\\
      0 &\frac{1}{4}(2H+i) e^{-v}
      \end{array}
\right).
$$

\noindent We denote by $M_1,M_2$ the matrices from the Gauss--Weingarten
equations \eqref{gauss_weingraten_eq}:  
\begin{equation*}
M_1=\begin{pmatrix}
     v_z&Be^{-v}\\
    -e^v&0
    \end{pmatrix},
\quad M_2=\begin{pmatrix}
    0 &e^v\\
    -\ovl{B}e^{-v}&v_{\ovl{z}}
    \end{pmatrix}.
\end{equation*}
%be the matrices of the 
%Gauss--Weingarten equations.
By differentiating the identity \eqref{v_cond}, we obtain that:
\begin{align}	
0&=\partial\big(\ovl{\psi}^t\wh\psi
           \big)
  =\ovl{\psi}^t\big(\ovl{M}_2^t\wh+\partial\wh+\wh M_1
               \big)\psi\nonumber\\
 &=\ovl{\psi}^tD_1\wh\psi
  =\tr\left[D_1\wh
            \left(\begin{array}{cc}
                  |\psi_1|^2&\psi_1\ovl{\psi_2}
                  \\
                  \ovl{\psi_1}\psi_2&|\psi_2|^2
                  \end{array}
            \right)
      \right],
\label{cond_1}\\                                  %(9)
0&=\ovp{\partial}\big(\ovl{\psi}^t\wh\psi
                 \big)=
   \ovl{\psi}^t\big(\ovl{M}_1^t\wh+\ovp{\partial}\wh+\wh M_2
               \big)\psi\nonumber\\
 &=\ovl{\psi}^tD_2\wh\psi=\tr
     \left[D_2\wh
           \left(\begin{array}{cc}
                 |\psi_1|^2&\psi_1\ovl{\psi_2}\\
                  \ovl{\psi_1}\psi_2&|\psi_2|^2
                 \end{array}
           \right)
     \right]
\label{cond_2}.                                    %(10)
\end{align}
The matrices $D_1\wh$ and $D_2\wh$ in
the identities \eqref{cond_1} and \eqref{cond_2}
have the form:
\begin{align}
D_1\wh &=\frac14e^{-v}\left[\begin{array}{cc}
                     0&\frac{i\tau B}
                            {|e^v|^2}\\
                     i\tau&\varkappa
                     \end{array}
               \right],\\                         %(11)
D_2\wh &=\frac14e^{-v}\left[\begin{array}{cc}
                     -\ovl{\varkappa}&i\sigma\\
                     \frac{i\sigma\ovl{B}}
                          {|e^v|^2}&0
                     \end{array}
               \right],                           %(12)
\end{align}
where coefficients 
$\tau$, $\sigma$ and $\varkappa$
can be expressed in terms of $v$ as follows: 
%that are given by the following equations
%$$
\begin{gather*}
i\tau=(2H-i)e^{\ovl{v}}-(2H+i)e^{v},\\
i\sigma=(2H-i)e^{v}-(2H+i)e^{\ovl{v}},\\
\varkappa=(2H+i)\big(\ovl{v}-v
                \big)_z .
\end{gather*}

%\begin{remark}                                 
%For minimal surfaces ($H=0$) the reality %condition is $\mathrm{Im}\,v=\pm\frac{\pi}2$ %and the matrices~$D_1\wh$ and $D_2\wh$ are %equal to zero, i.e., the %relations~\eqref{cond_1} and \eqref{cond_2} %are satisfied. If $H\ne 0$ then the imaginary %part of~$v$ is not constant and the
%coefficients of the matrices $D_1\wh$ and %$D_2\wh$ are not equal to zero.
%\end{remark}

\noindent Let $H=0$.  
Then the identity \eqref{cond1} implies 
that: 
\begin{equation} 
\label{real_cond_H_equals_0}   
   \operatorname{Re} (e^v)=0. 
\end{equation}	
If the identity \eqref{real_cond_H_equals_0} 
holds, then $\sigma=\tau=\varkappa=0$,  
so the matrices $D_1\wh$ and $D_2\wh$ are identically equal to zero. Therefore, 
there exists 
$\psi = \left(\begin{smallmatrix}
	\psi_1\\
	\psi_2
\end{smallmatrix}
\right)$ 
for which the Gauss--Weingarten 
equations and the identity:
$$
   e^v = \frac{i}{4} (|\psi_2|^2-|\psi_1|^2) 
$$  
are satisfied. 
For the case $H=0$ we refer to the equation 
\eqref{real_cond_H_equals_0} as the 
{\it reality condition}. 
%Note that the equations %\eqref{real_cond_H_equals_0} and 
%\eqref{compcond} are clearly 
%sufficient for the existence   
%of the generating spinor

%For minimal surfaces ($H=0$) the reality %condition is
%$\mathrm{Im}\,v=\pm\frac{\pi}2$ and the %matrices~$D_1\wh$ and $D_2\wh$ are equal to
%zero, i.e., the relations~\eqref{cond_1} and %\eqref{cond_2} are satisfied. If
%$H\ne 0$ then the imaginary part of~$v$ is not %constant and the
%coefficients of the matrices $D_1\wh$ and %$D_2\wh$ are not equal to zero

\noindent Let $H \ne 0$. 
%Without loss of generality we assume that 
We will require that 
$\psi_1\overline{\psi}_2\ne 0$ (equivalently, 
$\psi_1 \neq 0$ and $\psi_2 \neq 0$).
Recall that $\psi_1 \overline{\psi}_2$ is a complex coefficient 
in the decomposition \eqref{decomposition_f_z}: 
$f^{-1} \frac{\partial f}
{\partial z}=
\frac{i}{2}\big(\ovl{\psi}_2^2+\psi_1^2
\big)e_1+
\frac{1}{2}\big(\ovl{\psi}_2^2-\psi_1^2
\big)e_2+\psi_1\ovl{\psi_2} e_3$.     
We denote by $n$ a unit normal vector to a surface and by $E_3$ the left--invariant vector field   
generated by $e_3$. 
The assumption $\psi_1 \overline{\psi}_2 \neq 0$ holds if and only if: %at a given point 
$$n \neq \pm E_3.$$
It can be verified that this is equivalent 
to the following: 
\begin{equation} 
\label{first_assumption_v}   
   e^{\pm(\overline v  - v)} \neq \frac{2H-i}{2H+i}.	
\end{equation}	 
Note that there exists no surface in $\nil$ for which $n = \pm E_3$ (equivalently, $\psi_1\overline{\psi}_2 =0$) in some open neighborhood \cite{BT05}.

We denote by $\Omega$ the following matrix: 
$$
\Omega=
%\left(\begin{smallmatrix}
\left(\begin{array}{cc}
                 \,|\psi_1|^2&\psi_1\ovl{\psi_2}\\\noalign{\vskip2pt}
                 \ovl{\psi_1}\psi_2&\,|\psi_2|^2
%               \end{smallmatrix}
      \end{array}
\right).
 $$
Since $\psi_1 \overline{\psi}_2 \neq 0$, 
the matrix $\Omega$ can be represented as:
$$
\Omega=|\psi_1|^2 \left(\begin{array}{cc}
            1&\frac{\ovl{\psi_2}}
                   {\ovl{\psi_1}}\\
            \frac{\psi_2}
                 {\psi_1}&\,\,\big\vert\frac{\psi_2}
                                           {\psi_1}
                          \big\vert^2
      \end{array}
\right)= l\left(\begin{array}{cc}
        1&\ovl{\xi}\\
        \xi&|\xi|^2
       \end{array}
 \right),
$$
where $l=|\psi_1|^2$ and $\xi=\dfrac{\psi_2}{\psi_1}$. Then
the equations~\eqref{cond_1} and \eqref{cond_2} %can be represented as follows
take the following form:
%\vskip-7pt\noindent
\begin{align}	
&\frac{i\tau B}
      {|e^v|^2}\xi+i\tau\ovl{\xi}+\varkappa|\xi|^2=0,
      \label{equation1_xi}
\\
&\frac{i\sigma B}
      {|e^v|^2}\xi+i\sigma\ovl{\xi}+\varkappa=0.
\label{equation2_xi}
\end{align}
Note $\tau=0$ and $\sigma =0$ if and only if 
$e^{\overline v - v} =  \frac{2 H + i}{2 H - i}$
and $e^{\overline v - v} = \frac{2H-i}{2H+i}$, 
respectively. But in \eqref{first_assumption_v} 
we already have that $e^{\pm(\overline v  - v)} \neq \frac{2H-i}{2H+i}$, so $\tau \neq 0$ and 
$\sigma \neq 0$.  

We will additionally require that  
$\varkappa$ is 
not equal to zero (equivalently,  
$(\overline v - v)_z \ne 0$).  	
It can be verified that this is equivalent 
to the assumption that:
$$ 
  \frac{\partial}{\partial z} \langle n, E_3 \rangle  \neq  0. 
$$ 
Then the equations \eqref{equation1_xi} and 
\eqref{equation2_xi} imply that: 
\begin{equation}
\label{|xi|2_formula}
|\xi|^2=\frac{\tau}{\sigma}.
\end{equation}
The equations \eqref{equation1_xi}, \eqref{equation2_xi} and \eqref{|xi|2_formula} 
imply that: %\vskip-7pt\noindent
\begin{equation}
\label{equation_linear_xi}
\begin{split}
\frac{B}{|e^v|^2}\xi+\ovl{\xi}=- \frac{\varkappa}
     {i\sigma},
\\
\xi+\frac{\ovl{B}}
         {|e^v|^2}\ovl{\xi}=
    \frac{\ovl{\varkappa}}
         {i\sigma}.
\end{split} 
\end{equation}
%We now solve this system with respect to 
Resolving the system \eqref{equation_linear_xi} with respect to $\xi$ we obtain: 
\begin{equation}
\label{formula_for_xi}	
\xi=\frac{\ovl{\varkappa}+
	\varkappa\frac{\ovl{B}}
         {|e^v|^2}}
         {i\sigma\left(1-\frac{|B|^2}
                              {|e^v|^4}
                 \right)}.
\end{equation}
The equations \eqref{|xi|2_formula} and
\eqref{formula_for_xi} imply that:
\begin{equation}\label{real_cond1}                %(13)
\left(\ovl{\varkappa}+\varkappa
      \frac{\ovl{B}}
           {|e^v|^2}
\right) \left(\varkappa+\ovl{\varkappa}
      \frac B{|e^v|^2}
\right)=\tau\sigma \left(1-\frac{|B|^2}
             {|e^v|^4}
\right)^2 .
\end{equation}

\noindent We say that a given point of a surface in $\nil$ is {\it umbilic} if at this point $\wta  = 0$. We consider constant 
mean curvature surfaces in neighborhoods
of non--umbilic points ($\wta \ne 0$). 
In some neighborhood of 
a non--umbilic point of a constant mean 
curvature surface we can always choose a conformal parameter $z$ such that:
\begin{equation} 
\label{choice_of_B}  
  B=\dfrac{2H+i}{2H-i}.
\end{equation}          
%\begin{remark}                                 
%We consider a surface only in a neighborhood %of a non--umbilical point,
%i.~e., under the condition  $\wta\ne 0$. The %surfaces with vanishing
%quadratic differential $\wta dz^2$ are known. %For instance, if
%$H\ne 0$ then these are the spheres of %constant mean curvature~\cite{AR05}.
%\end{remark}
%By changing a conformal parameter we may put
%that:
%Locally we can choose a conformal parameter 
%$z$ such that:  
%$$B=\dfrac{2H+i}{2H-i}.$$ 
Assuming \eqref{choice_of_B}, the identity \eqref{real_cond1} 
takes the form: 
%is equivalent to the following equation:
%$$
\begin{equation}\label{real_cond2}              %(14)
%\hskip-2mm
\begin{split}
\left\vert(\ovl{v}{-}v)_z{-}(\ovl{v}{-}v)_{\ovl{z}}     \frac1{|e^v|^2}
\right\vert^2=
\,\,\,\,\,\,\,\,\,\,\,\,\,\,\,\,\,\,\,\,
\,\,\,\,\,\,\,\,\,\,\,\,\,\,\,\,\,\,\,\,
\,\,\,\,\,\,\,\,\,\,\,\,\,\, 
\\
\left(e^{2v}{+}e^{2\ovl{v}}{-}2\,
\operatorname{Re}
\left(\frac{2H+i}{2H-i}\right)|e^v|^2
\right)   
\left(1{-}\frac1{|e^v|^4} \right)^2.
\end{split}
\end{equation}
%Let $v=\rho+i\varphi$. 
Let $\rho = \operatorname{Re}(v)$ and 
$\varphi = \operatorname{Im}(v)$. 
Then $\ovl{v}-v=-2i\varphi$ and
$|e^v|=e^{\rho}$. 
%In this case 
In terms of $\rho$ and $\varphi$ 
the identity \eqref{real_cond2} 
finally takes the form:
\begin{equation}
\label{real_cond3} 
\frac{\varphi_x^2}
     {\cosh^2\rho}+
\frac{\varphi_y^2}
     {\sinh^2\rho}=8
\left(\cos 2\varphi-
\operatorname{Re}\left(\frac{2H+i}{2H-i}
\right)\right),
\end{equation}
where $z= x+iy$. Similarly to minimal surfaces ($H = 0$),  for the case $H \ne 0$ 
we will refer to the equation 
\eqref{real_cond3} as the reality 
condition.  
%Finally we formulate the following theorem.
We summarize the obtained results in Theorems \ref{theocmcsyst} and 
\ref{minimal_surf_thm}.  

\begin{theorem}\label{theocmcsyst}                %Теорема 1
A surface of nonzero constant mean curvature 
in $\nil$ in some neighborhood of
a non--umbilic point, 
where we assume that 
a unit normal vector $n \neq E_3$
%$e^{\pm 2 i \varphi} \ne 
%\frac{2H+i}{2H-i}$ 
and 
%$\varphi_z \ne 0$, 
$\frac{\partial}{\partial z} \langle n, 
E_3 \rangle \ne 0$,  
corresponds to a solution $v=\rho+i\varphi$
of the following system:
$$
\begin{cases}
%v_{z\ovl{z}}+2\sinh 2v=0,
\Delta v+ 8 \sinh 2v=0,
\\
\dfrac{\varphi_x^2}
     {\cosh^2\rho}+
\dfrac{\varphi_y^2}
     {\sinh^2\rho}=8
\left(\cos 2\varphi-\operatorname{Re}
\left(\dfrac{2H+i}{2H-i}\right) \right).
\end{cases}
$$
\end{theorem}

\begin{theorem} 
\label{minimal_surf_thm}	
   A minimal surface in $\nil$ 
   in some neighborhood 
   of a non--umbilic point corresponds to a 
   solution of the following system:
   $$
   \begin{cases}
   \Delta v+ 8 \sinh 2v=0,
   \\ 
   \operatorname{Re}(e^v)=0.
   \end{cases}
   $$ 
\end{theorem}	

\section*{Acknowledgments}

The author thanks Iskander Taimanov and Uwe Abresch for useful conversations.

\end{document}